\documentclass[letterpaper, 10 pt, conference]{ieeeconf}  % Comment this line out if you need a4paper

\IEEEoverridecommandlockouts                              % This command is only needed if 
\usepackage{cite}
\usepackage{amsmath,amssymb,amsfonts}
  \usepackage{amsthm}
\usepackage{algorithmic}
\usepackage{graphicx}
\usepackage{textcomp}
\usepackage{xcolor}
\usepackage{mathrsfs}
\usepackage{array}
\usepackage{multirow} 
\usepackage{graphicx}
\usepackage{algorithm}
\usepackage{algorithmic}
\usepackage{float}
\usepackage{color,xcolor}
\usepackage{array}
\usepackage{arydshln}
\usepackage{multirow}

\def \R {\mathbb{R}}
\def \B {\mathcal{B}}

\def \N {\mathcal{N}}
\def \bx {{\bf x}}

\def \bk {{\bf k}}

\def \Z {\mathcal{Z}}
\def \P {\mathbb{P}}

\theoremstyle{plain}
\newtheorem{THEOREM}{Theorem}[section]

\newtheorem{theorem}[THEOREM]{Theorem}
\newtheorem{proposition}[THEOREM]{Proposition}
\newtheorem{corollary}[THEOREM]{Corollary}

\theoremstyle{remark}
\newtheorem{example}[THEOREM]{Example}

\allowdisplaybreaks[4]

\title{\LARGE \bf
A Decomposition Approach for the Gain Function in the Feedback Particle Filter*
}

\author{Ruoyu Wang$^{1}$, Huimin Miao$^{2}$ and Xue Luo$^{3,\sharp}$, {\it IEEE senior member}% <-this % stops a space
\thanks{*This work is financially supported by the National Key R\&D Program of China (Grant No. 2022YFA1005103), National Natural Science Foundation of China (Grant No. 12271019), and Natural Science Foundation of Henan Province (Grant No. 252300420900).}% <-this % stops a space
\thanks{$^{1}$R. Wang is with School of Mathematical Sciences, Beihang University, Beijing, P. R. China 102206. 
        {\tt\small WangRY@buaa.edu.cn}}%
\thanks{$^{2}$H. Miao is with  School of Mathematics and Information Science, Henan Polytechnic University, Jiaozuo, Henan, P. R. China 454003.
        {\tt\small 13269418722@163.com}}%
\thanks{$^{3}$X. Luo is with School of Mathematical Sciences and Key Laboratory of Mathematics, Informatics and Behavioral Semantics (LMIB), Beihang University, Beijing, P. R. China 100191.
        {\tt\small xluo@buaa.edu.cn}}%
\thanks{$^{\sharp}$X. Luo is the corresponding author.}
}

\begin{document}

\maketitle
\thispagestyle{empty}
\pagestyle{empty}

%%%%%%%%%%%%%%%%%%%%%%%%%%%%%%%%%%%%%%%%%%%%%%%%%%%%%%%%%%%%%%%%%%%%%%%%%%%%%%%%
\begin{abstract}
    The feedback particle filter (FPF) is an innovative, control-oriented and resampling-free adaptation of the traditional particle filter (PF). In the FPF, individual particles are regulated via a feedback gain, and the corresponding gain function serves as the solution to the Poisson's equation equipped with a probability-weighted Laplacian. Owing to the fact that closed-form expressions can only be computed under specific circumstances, approximate solutions are typically indispensable. This paper is centered around the development of a novel algorithm for approximating the gain function in the FPF. The fundamental concept lies in decomposing the Poisson's equation into two equations that can be precisely solved, provided that the observation function is a polynomial. A free parameter is astutely incorporated to guarantee exact solvability. The computational complexity of the proposed decomposition method shows a linear correlation with the number of particles and the polynomial degree of the observation function. We perform comprehensive numerical comparisons between our method, the PF, and the FPF using the constant-gain approximation and the kernel-based approach. Our decomposition method outperforms the PF and the FPF with constant-gain approximation in terms of accuracy. Additionally, it has the shortest CPU time among all the compared methods with comparable performance.
\end{abstract}

%%%%%%%%%%%%%%%%%%%%%%%%%%%%%%%%%%%%%%%%%%%%%%%%%%%%%%%%%%%%%%%%%%%%%%%%%%%%%%%%
\section{Introduction}\label{sec-1}

Nonlinear filtering (NLF) is a sophisticated field that focuses on extracting valuable information from sensor-corrupted data. It forms the bedrock of a wide array of applications, such as target tracking and navigation, air traffic management, and meteorological monitoring, etc.

The NLF problem for diffusion processes is formulated in terms of the following It\^o stochastic differential equation (SDE):
\begin{equation}\label{eqn-NLF}
\left\{\begin{aligned}
dX_t=&g(X_t)dt+\sigma(X_t)dB_t\\
dZ_t=&h(X_t)dt+dW_t
\end{aligned}\right..
\end{equation}
Here, $X_t\in\R^d$ denotes the state, and $Z_t\in\R^m$ represents the observation. $\{B_t\},\{W_t\}$ are two mutually-independent standard Wiener processes taking values in $\R^d$ and $\R^m$ respectively. The functions $g(\cdot):\ \R^d\rightarrow\R^d$, $h(\cdot):\ \R^d\rightarrow\R^m$, and $\sigma(\cdot):\ \R^d\rightarrow\R^{d\times d}$ are bounded continuous. Additionally, $g$ and $\sigma$ satisfy the Lipschitz condition. The overarching objective of the NLF problem is to approximate the posterior distribution $p_t^*(\bx)=\P(X_t|\Z_t)$, where $\Z_t=\sigma\{Z_s:\,s\leq t\}$ is the $\sigma$-field generated by the observation history up to time $t$.

The investigation of the NLF algorithm dates back to the 1960s, marking a long and rich history of research in this domain. In recent years, Yang et al. \cite{YMM:13} introduced a control-oriented, resampling-free algorithm known as the feedback particle filter (FPF). This innovative algorithm artfully combines the feedback structure, characteristic of the Kalman filter, with the flexibility offered by the ensembles in the PF. In the FPF framework, each individual particle is governed by a controlled SDE, expressed as:
\begin{equation}\label{eqn-2.20}
dX_t^i = g(X_t^i)dt+\sigma(X_t^i)dB_t^i+K(X_t^i,t)dZ_t+u(X_t^i,t)dt.
\end{equation}
Here, the gain function $K:\R^d\times\R_+\rightarrow\R^{d\times m}$ plays a crucial role. It serves as the solution to a Poisson's equation given by:
\begin{equation}\label{eqn-BVP}
	\nabla^T(p_t^*K)=-\left(h-\hat h\right)^Tp_t^*,
\end{equation}
where $\hat h=\int_{\R^d}h(\bx)p_t^*(\bx)d\bx$. The symbol $\nabla^T$ represents the divergence operator. Additionally, the function $u$ in the SDE \eqref{eqn-2.20} is defined as:
\begin{equation}\label{eqn-u}
u=-\frac12K\left(h+\hat h\right)+\frac12\sum_{k=1}^d\sum_{s=1}^mK_{ks}\frac{\partial K_{ls}}{\partial x_k}.
\end{equation}
The consistency of the FPF has been rigorously established in \cite{YMM:13}. This theoretical result states that if the initial particles $\{X_0^i\}_{i=1}^{N_p}$ are sampled from the initial density $p_0^*(\bx)$ of $X_0$, then for any Borel measurable set $A\in\B(\R^d)$, the empirical distribution of the controlled particles, defined as
\[
    p_t^{N_p}(A):=\frac1{N_p}\sum_{i=1}^{N_p} 1_{\{X_t^i\in A\}},
\]
approximates the true posterior density $p_t^*(A)$ for each $t$, as the number of the particles $N_p$ approaches infinity. 

From a numerical perspective, Surace et al. \cite{SKP:19} has demonstrated through experiments that when the gain function is accurately determined, the performance of the FPF can surpass that of almost all traditional NLF algorithms. However, as pointed out in \cite{TM:16}, finding the solution to \eqref{eqn-BVP} remains a formidable challenge. This is because, in general, there is no closed-form solution available for it. Furthermore, the true density $p_t^*(\bx)$ is not explicitly known, adding to the complexity of the problem. Consequently, there is a pressing need to develop an accurate and efficient method for obtaining an approximate gain function.

In the quest for this issue, several methods have been developed. Yang et al. \cite{YLMM:13} first proposed the constant-gain approximation, which is a basic average of particle controls. The Galerkin approximation, developed in \cite{YMM:13, YLMM:16}, improves upon it but requires pre-defined basis functions, scaling poorly with high-dimensional states and potentially suffering from the Gibbs phenomena as noted in \cite{TM:16}. To address the selection of basis functions, Berntorp et al. \cite{BG:16} introduced the proper orthogonal decomposition-based (POD-based) method for adaptive selection. Taghvaei et al. \cite{TM:16} then proposed the kernel-based algorithm, eliminating the need for basis functions entirely. Error analysis of this kernel-based method was explored in \cite{TMM:17,TMM:20}. Berntorp \cite{B:18} compared the constant-gain, POD-based, and kernel-based approximations in the FPF, offering insights into their performances.

In this paper, we introduce a novel decomposition method aimed at obtaining an approximate solution to \eqref{eqn-BVP} subject to a specific boundary condition. This boundary condition is different from the zero-mean constraint in \cite{TM:16}. Without loss of generality, we make the assumption:
\begin{quote}
    (As) \quad the observation $h(\bx)$ is a polynomial in $\bx$.
\end{quote}
In cases where $h(\bx)$ does not satisfy this condition, we can initially approximate it using a polynomial function. 

The fundamental concept of our proposed decomposition method is structured into the following two steps:
\begin{enumerate}
    \item[(1)] The density function $p_t^*$ is approximated by
    \begin{equation}\label{eqn-4}
    p_t^{N_p,\Sigma}(\bx)=\frac1{N_p}\sum_{i=1}^{N_p}\N(\bx;X^i_t,\Sigma),
    \end{equation}
    where \begin{align}\label{eqn-Gaussian multi}
	&\mathcal{N}(\bx;X_t^i,\Sigma)\\\notag
    =&\frac{1}{(2\pi)^{\frac{d}{2}}|\Sigma|^{\frac{1}{2}}}\exp\left\{-\frac{1}{2}(\bx-X_t^i)^T\Sigma^{-1}(\bx-X_t^i)\right\}
  \end{align}
  represents the multivariate Gaussian density function.
    \item[(2)] The equation
    \begin{equation}\label{eqn-62}
    \nabla^T\left(p_t^{N_p,\Sigma}K\right)=-\left(h-\hat h^{N_p}\right)^Tp_t^{N_p,\Sigma}
    \end{equation}
  is decomposed into two components that are exactly solvable, where
   \begin{equation}\label{eqn-hat hN}
    \hat h^{N_p}=\int_{\R^d}h(\bx)p^{N_p,\Sigma}_t(\bx)d\bx.
    \end{equation}
     Given appropriate boundary condition, explicit solutions for these components can be derived, see Section \ref{sec-observation}.
\end{enumerate}

It is important to emphasize that in this decomposition method, the only source of approximation lies in Step (1) described above, while Step (2) is solved exactly. We shall explain the reasonableness of the approximation in Step (1) in section \ref{sec-2}.

The organization of this paper is as follows: In Section \ref{sec-2}, we review some preliminary concepts. Section \ref{sec-3} presents the decomposition method to obtain the gain function in the FPF, and Section \ref{sec-4} details the numerical experiments to verify our algorithm.

\section{Preliminary}\label{sec-2}

The accuracy of the approximate gain function $K$ plays a pivotal role in determining the performance of the FPF. The gain function $K$ is governed by Poisson's equation \eqref{eqn-BVP}. When approaching the analysis of this equation, a fundamental consideration is its well-posedness, which encompasses the existence and uniqueness of the solution. This aspect is not only a theoretical cornerstone but also has practical implications for the reliability and effectiveness of the FPF. Ensuring well-posedness provides a solid foundation for accurately approximating the gain function and subsequently optimizing the performance of the overall filtering process.

Laugesen et al. \cite{LMMR:15} demonstrated the existence and uniqueness of the solution to \eqref{eqn-BVP} under the substitution $K=\nabla\phi$ and the zero-mean constraint $\int_{\R^d}p_t^*(\bx)\phi(\bx) d\bx=0$. This result holds given specific conditions on the density $p_t^*$ and the smoothness of the observation function $h$. However, it should be noted that the conditions imposed on $p_t^*$ are primarily for the sake of rigorous mathematical proof. In practice, $p_t^*$ is entirely unknown rendering these conditions unverifiable. In this paper, rather than the zero-mean constraint, we impose an appropriate boundary condition and introduce a substitution for $p_t^*$.

For the boundary condition, we perform an integration of \eqref{eqn-BVP} over $\R^d$ with respect to $\bx$:
\begin{equation}\label{eqn-61}
    \int_{\R^d}\nabla^T(p^*_tK)d\bx=-\int_{\R^d}(h-\hat h)^Tp_t^*d\bx=0.
\end{equation}
Invoking the Divergence theorem, it becomes evident that it is reasonable to impose the following constraint:
\begin{equation}\label{eqn-6}  \lim_{R\rightarrow\infty}\int_{S^d(R)}p^*_tK\cdot d\sigma
    =\lim_{R\rightarrow\infty}\int_{B^d(R)}\nabla^T(p^*_tK)dx=0,
\end{equation}
where $B^d(R)=\{\bx\in\R^d:|\bx|\leq R\}$, and $S^d(R)$ represents the surface of the ball $B^d(R)$. In the scalar scenario, when $d=1$, \eqref{eqn-6} reduces to the boundary condition $p_t^*K(-\infty)+p_t^*K(+\infty)=0$. Consequently, we may impose the boundary condition 
\begin{equation}\label{eqn-41}
 p_t^*K(-\infty)=p_t^*K(+\infty)=0.
\end{equation}

To render \eqref{eqn-BVP} tractable, it is necessary to seek an approximation of $p_t^*$. The approximation $p_t^{N_p,\Sigma}(\bx)$ given by \eqref{eqn-4} is a reasonable choice. In fact, the multivariate Gaussian density $\N(\bx;X_t^i,\Sigma)$ is commonly employed to approximate the Dirac delta function $\delta_{X_t^i}(\bx)$. In the sense that as the covariance matrix $\Sigma$ degenerates to zero-matrix,  $\N(\bx;X_t^i,\Sigma)$ converges to $\delta_{X_t^i}(\bx)$ in the weak sense, as detailed in Chapter 2, \cite{GS:68}. This weak-convergence property provides a theoretical foundation for using $p_t^{N_p,\Sigma}(\bx)$ as an approximation for $p_t^*(\bx)$.

Consequently, in Section \ref{sec-3}, our efforts will be concentrated on devising an exact solution for \eqref{eqn-62} that has an explicit form in the scalar case. Subsequently, we shall proceed to verify the condition stated in \eqref{eqn-41}. 

\section{The Decomposition Method}\label{sec-3}

As described in \cite{YLMM:13}, consider a vector-valued function $\varphi(\bx)=[\varphi_{1}(\bx),\cdots,\varphi_m(\bx)]$ with $m$ components, such that for each $j=1,\cdots,m$, the $j$-th column of the matrix $K$, denoted as $K_{\cdot j}(\bx)$, satisfies $K_{\cdot j}=\nabla\varphi_j(\bx)$. Under this definition, \eqref{eqn-62} is then equivalent to  
\begin{align}\label{eqn-pN_phi}
	\nabla^T\left(p_t^{N_p,\Sigma}(\bx)\nabla\varphi_j(\bx)\right)=-\left(h_j(\bx)-\hat h_j^{N_p}\right)p_t^{N_p,\Sigma}(\bx),
\end{align}
for $j=1,\cdots,m$. In this section, our objective is to develop a decomposition method that divides \eqref{eqn-pN_phi} into two parts, both of which admit explicit solutions. This decomposition approach represents, under the assumption (As), a highly efficient means of obtaining the {\it exact} solution. Let us assume that for all $j=1,\cdots,m$, the function $h_j(\bx)$ is a polynomial of degree less or equal to $p$. We can then express $h_j(\bx)$ as
	\begin{equation}\label{eqn-3.14}
	h_j(\bx)=\sum_{\bk\in\Omega}a_{j,\bk}H_{\bk}(\bx),
	\end{equation}
for each $j=1,\cdots,m$. Here, $\Omega:=\{\bk\in\mathbb{N}^d:\,0\leq|\bk|_1\leq p\}$, with the $l_1$-norm of $\bk$ defined as $|\bk|_1=\sum_{l=1}^dk_l$, and $H_{\bk}(\bx)=\Pi_{l=1}^dH_{k_l}(\bx)$ represents the multi-variable Hermite polynomials.

\subsection{The decomposition of \eqref{eqn-62}}\label{sec-3.1}

In this subsection, we shall first decompose \eqref{eqn-62} into two components. Each of these components can be solved both precisely and with high computational efficiency, thus enabling us to make substantial headway in addressing the equation within the context of our overall problem-solving framework. This decomposition step is fundamental to our approach, as it allows us to break down a complex equation into more manageable sub-problems, for which exact solutions can be derived, under the assumption (As).

\begin{proposition}\label{prop-1}For each $j = 1,\cdots,m$, the gain function $\nabla\varphi_j(\mathbf{x}) $ in \eqref{eqn-pN_phi} is given by
\begin{align}\label{eqn-3.10}\notag
    \nabla\varphi_j(\mathbf{x})=&\frac{1}{\sum_{i = 1}^{N_p}\mathcal{N}(\mathbf{x};X_t^i,\Sigma)}\\
    &\cdot\sum_{i = 1}^{N_p}\left[\nabla_{\vec{\lambda}}\psi_j^i(\mathbf{x})+\mathcal{N}(\mathbf{x};X_t^i,\Sigma)\nabla\varphi_j^i(\mathbf{x})\right],
\end{align}
where $\vec{\lambda}=(\lambda_1,\cdots,\lambda_d) $ represents the eigenvalues of $\Sigma^{-1}$, and $\nabla_{\vec{\lambda}}\psi_j^i(\mathbf{x}):=\left(\frac{1}{\lambda_l}\frac{\partial\psi_j^i(\mathbf{x})}{\partial x_l}\right)_{l = 1}^{d}$. The functions $\varphi_j^i(\mathbf{x}) $ and $\psi_j^i(\mathbf{x}) $ satisfy the following equations:
\begin{equation}\label{eqn-3.1}
    \nabla^T\left(\mathcal{N}(\mathbf{x};X_t^i,\Sigma)\nabla\varphi_j^i(\mathbf{x})\right)
    =\left(-h_j(\mathbf{x})+C_j^i\right)\mathcal{N}(\mathbf{x};X_t^i,\Sigma),
\end{equation}
for any constant $C_j^i$, and \begin{equation}\label{eqn-3.2}
    \nabla^T\left(\nabla_{\vec{\lambda}}\psi_j^i(\mathbf{x})\right)
    =\left(\hat{h}_j - C_j^i\right)\mathcal{N}(\mathbf{x};X_t^i,\Sigma),
\end{equation}
respectively.
\end{proposition}

\begin{proof}
   Given $j=1,\cdots,m$, we have
   \begin{align*}
    &\nabla^T\left[\frac{1}{N_p}\sum_{i=1}^{N_p}\mathcal{N}(\bx;X_t^i,\Sigma)\nabla\varphi_j(\bx)\right]\\
    \overset{\eqref{eqn-4},\eqref{eqn-pN_phi}}=&\frac1{N_p}\sum_{i=1}^{N_p}-\left(h_j(\bx)-\hat h_j^{N_p}\right)\N(\bx;X_t^i,\Sigma)\\
    \overset{\eqref{eqn-3.1},\eqref{eqn-3.2}}=&\frac{1}{N_p}\sum_{i=1}^{N_p}\nabla^T\left[\mathcal{N}(\bx;X_t^i,\Sigma)\nabla\varphi^i_j(\bx)\right]\\
    &+\frac{1}{N_p}\sum_{i=1}^{N_p}\nabla^T\left(\nabla_{\vec\lambda}\psi_j^i(\bx)\right)\\
    =&\nabla^T\left[\frac{1}{N_p}\sum_{i=1}^{N_p}\left(\mathcal{N}(\bx;X_t^i,\Sigma)\nabla\varphi^i_j(\bx)+\nabla_{\vec\lambda}\psi_j^i(\bx)\right)\right].
\end{align*}
Consequently, \eqref{eqn-3.10} follows immediately.
\end{proof} 

According to Proposition \ref{prop-1}, once obtaining the {\it exact} solutions to \eqref{eqn-3.2} and \eqref{eqn-3.1}, an explicit formula for the gain function is provided by \eqref{eqn-3.10}. However, the selection of the constant $C_j^i$ is a crucial aspect that requires careful consideration. Our proposed strategy for determining $C_j^i$ consists of the following two steps: 
\begin{itemize}
	\item[1)] We aim to select the constant $C_j^i$ in such a way that \eqref{eqn-3.1} can be precisely solved using the Galerkin spectral method via Hermite polynomials. We use the backward recursion to solve for the Hermite coefficients, under the assumption (As). 
	\item[2)] After choosing the appropriate $C_j^i$ in Step 1), we then seek to obtain an exact solution for \eqref{eqn-3.2} with the same value of $C_j^i$. 
\end{itemize}

It should be noted that the Galerkin spectral method had been previously employed to solve \eqref{eqn-BVP} in the context where $p_t^*$ was approximated by the empirical density of the particles within the zero-mean Sobolev space, as reported in \cite{TM:16}. However, our approach differs from theirs, we introduce a free parameter $C_j^i$ that can be judiciously selected. This selection of $C_j^i$ enables us to render \eqref{eqn-3.1} exactly solvable within the space spanned by the Hermite polynomials. This unique feature of our method is significant in the sense that it is a possibility of completely surmounting the curse of dimensionality, as discussed in \cite{SKP:19}. 

As the initial phase of our decomposition method, we shall test it by seeking exact solutions to \eqref{eqn-3.1} and \eqref{eqn-3.2} in the one-dimensional ($d = 1$) scenario. This allows us to better understand how this method work. The one-dimensional case serves as a crucial starting point before applying the method to higher-dimensional problems. The decomposition method in high-dimensional scenario will be discussed briefly in Section \ref{sec-3.c}.

\subsection{The decomposition method in the scalar case}\label{sec-observation}

It is evident that when the dimension $d = 1$, as stated in \cite{YMM:13}, \eqref{eqn-pN_phi} can be solved through direct integration. At first glance, this might make our proposed decomposition method seem unnecessary. However, it is important to note that for dimensions $d\geq2$, the presence of the divergence operator $\nabla^T$ makes direct integration inapplicable.

At the moment, let us consider the case where the dimension $d = 1$. In accordance with the assumption (As), we have the polynomial representation of the function $h(x)$ given by:
\begin{equation}\label{eqn-42}
    h(x)=\sum_{k=0}^pa_kH_{k}(x),
\end{equation}
    for some $p\geq1$. In the scalar scenario, instead of solving for the function $\varphi$ as in Proposition \ref{prop-1}, we shall directly find the gain function $K$. It is straightforward to observe that \eqref{eqn-pN_phi} simplifies to:
\begin{equation}\label{eqn-70}
	\left(p_t^{N_p,\epsilon}(x)K(x)\right)'=-\left(h(x)-\hat h^{N_p}\right)p_t^{N_p,\epsilon}(x),
\end{equation} 
where $p_t^{N_p,\epsilon}(x)=\frac1{N_p}\sum_{i=1}^{N_p}\N(x;X_t^i,\epsilon)$. The boundary condition \eqref{eqn-41} is adjusted such that $p_t^*$ is replaced by $p^{N_p,\epsilon}$. As a result of these simplifications in the scalar case, Proposition \ref{prop-1} directly leads to the following corollary.
\begin{corollary}
The gain function $K(x)$ in \eqref{eqn-70} is given by
\begin{align}\label{eqn-50}
	K(x)
    =&\frac1{\sum_{i=1}^{N_p}\mathcal{N}(x;X_t^i,\epsilon)}\sum_{i=1}^{N_p}\left[\mathcal{N}(x;X_t^i,\epsilon)K^i(x)+K^i_0(x)\right],
\end{align}
where the functions $K^i(x)$ and $K_0^i(x)$  satisfy
\begin{equation}\label{eqn-1}
	\left[\N(x;X_t^i,\epsilon){K^i}(x)\right]'=-\left(h(x)-C^i\right)\N(x;X_t^i,\epsilon),
\end{equation}
 and 
\begin{equation}\label{eqn-71}
	{K_0^i}'(x)=\left(\hat h^{N_p}-C^i\right)\N(x;X_t^i,\epsilon),
\end{equation}
for any constant $C^i$, respectively.
\end{corollary}

As elaborated in Section \ref{sec-3.1}, our approach involves a two-step process. First, we utilize the Hermite spectral method to solve \eqref{eqn-1} for $K^i$ with carefully chosen constant $C^i$. Subsequently, we proceed to solve for $K_0^i$ in \eqref{eqn-71} with the same value of $C^i$. It is of utmost importance to note that we choose $C^i$ precisely so that \eqref{eqn-1} can be solved exactly within the space spanned by the Hermite polynomials.

\subsubsection{For \eqref{eqn-1}} Employing the Hermite spectral method, we search for the solution $K^i(x)$ within the space that is spanned by Hermite polynomials of degree up to $q\geq p$, where $p$ is the highest polynomial degree of $h(x)$. We can represent $K^i(x)$ in the following form:
\begin{equation}\label{eqn-52}
	K^i(x)=\sum_{l=0}^q\hat K^i_lH_l(x),
\end{equation}
 where $H_l(x)$ denotes the Hermite polynomial of degree $l$, and $\hat K^i_l$ represents the corresponding coefficient associated with that polynomial. Through direct computations, we find that \eqref{eqn-1} is equivalent to:
\begin{align}\label{eqn-2}
	{K^i}'(x)-\frac{x-X_t^i}{\epsilon}K^i(x)\overset{\eqref{eqn-42}}=&-\sum_{k=0}^pa_kH_{k}(x)+C^i.
\end{align}
When we substitute \eqref{eqn-52} into \eqref{eqn-2}, we obtain that
\begin{align}\label{eqn-3.11}\notag
	&\sum_{l=0}^q\hat K_l^i\left(2-\frac 1{\epsilon}\right)lH_{l-1}(x)+\frac{X_t^i}{\epsilon}\sum_{l=0}^q\hat K^i_lH_l(x)\\
        &-\frac1{2\epsilon}\sum_{l=0}^q\hat K_l^iH_{l+1}(x)
        =-\sum_{l=0}^qa_lH_l(x)+C^i.
\end{align}
The derivation of the first two terms in \eqref{eqn-3.11} is based on the three-term recurrence property of Hermite polynomials. That is,
\begin{align*}
    {H}_{0}(x)=&1,\quad {H}_{1}(x)=2x,\\
{H}_{n+1}(x)=&2x{H}_{n}(x)-2n{H}_{n-1}(x),
\end{align*}
and for $n\geq1$, ${H}_{n}'(x)=2n{H}_{n-1}(x)$.
   
\begin{figure}[h!]
    \centering
    \includegraphics[trim = 0mm 28mm 0mm 16mm, clip,scale=1]{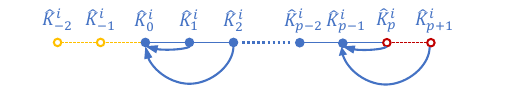}
    \caption{The coefficients $\hat K_k^i$ in \eqref{eqn-92}, $k=p-1,\cdots,1,0$, are calculated backwardly with the initial values $\hat K_p^i=\hat K_{p+1}^i\equiv0$ and the convention that $\hat K_k^i\equiv0$ when $k<0$. %The circles mark $\hat K_k^i=0$, while the dots represent $\hat K_k^i\not=0$ for $k\in\mathbb{Z}$.
    }\label{fig-d1}
\end{figure}

Multiply both sides of \eqref{eqn-3.11} by $H_k(x)e^{-x^2}$ for $k=0,1,2,\cdots$, and then integrate over $\R$ with respect to $x$. Leveraging the orthogonality property of Hermite polynomials, given by $\displaystyle\int_{-\infty}^{+\infty}H_m(x)H_n(x)\omega(x)dx=\gamma_n\delta_{mn}$, where $\gamma_n=\sqrt{\pi}2^nn!$ and
$
\delta_{mn}=\left\{\begin{aligned}
1,\quad m=n\\
0,\quad m\neq n
\end{aligned}\right.,
$
we obtain the following recurrence relation:
\begin{align}\label{eqn-54}\notag
	&\hat K_{k+1}^i\left(2-\frac1{\epsilon}\right)(k+1)+\frac{X_t^i}{\epsilon}\hat K_k^i-\frac1{2\epsilon}\hat K_{k-1}^i\\
        =&-a_k+C^i\delta_{k0},
\end{align}
for $k=0,1,2,\cdots$. Here, we adopt the convention that $\hat K_k^i\equiv0$, when $k<0$. It is important to note that the coefficients $a_k$s are non-zero only up to degree $p$. We can find an explicit solution for $\hat K_k^i$ with $k=0,1,\cdots,p-1$, by using backward recursion. We set the initial conditions $\hat K_k^i\equiv0$ for all $k\geq p$. The coefficient $\hat K_k^i$ is uniquely determined in a backward recursive manner by $\hat K_{k+1}^i$ and $\hat K_{k+2}^i$ through the formula:
\begin{equation}\label{eqn-92}
	\hat K_k^i=2\epsilon a_{k+1}+2(\epsilon-1)(k+2)\hat K_{k+2}^i+2X_t^i\hat K_{k+1}^i,
\end{equation}
for $k=p-1,\cdots,1,0$, with the initial values $\hat K^i_{p+1}=\hat K^i_{p}\equiv0$. Finally, we can determine the constant $C^i$ using $\hat K_0^i$ and $\hat K_1^i$ according to the formula:
\begin{equation}\label{eqn-Ci}
	C^i:=a_0+\frac{X_t^i}{\epsilon}\hat K^i_0+\left(2-\frac1{\epsilon}\right)\hat K^i_1.
\end{equation}
For the sake of clarity, the backward recursion process of the coefficients is visually depicted in Fig. \ref{fig-d1}.  

\subsubsection{For \eqref{eqn-71}}\label{sec-3.2.2} With the constant $C^i$ obtained in \eqref{eqn-Ci}, the solution can be derived through direct integration:
\begin{align}\label{eqn-59}\notag
	K_0^i(x)=&K_0^i(-\infty)+\int_{-\infty}^x\left(\hat h^{N_p}-C^i\right)\mathcal{N}(y;X_t^i,\epsilon)dy\\
	=&K_0^i(-\infty)+\left(\hat h^{N_p}-C^i\right)\frac12\left[1+\textup{erf}\left(\frac{x-X_t^i}{\sqrt{2\epsilon}}\right)\right].
\end{align}
Here, $\textup{erf}(x):=\frac 2{\sqrt\pi}\int_0^xe^{-\eta^2}d\eta$ represents the error function. In the context of the multivariate scenario, the task of determining the exact solution for \eqref{eqn-3.2}, which is the corresponding part of \eqref{eqn-71}, demands further exploration. This is the key part to make the decomposition algorithm tractable in high-dimensional NLF problems. 

\subsubsection{Boundary condition \eqref{eqn-41}} The first term in \eqref{eqn-50} that includes $K^i(x)$ naturally satisfies the boundary condition. This is due to the fact that the Gaussian function $\N(x;X_t^i,\epsilon)$ decays exponentially as $|x|$ approaches infinity, and $K^i(x)$ is a polynomial of degree at most $p$ for each $i=1,\cdots,N_p$. The boundary condition \eqref{eqn-41} allows us to analyze the behavior of the second term in \eqref{eqn-50} involving $K_0^i(x)$. First, we find that
\begin{align}\label{eqn-73}
	0\overset{\eqref{eqn-41}}=p_t^{N_p,\epsilon}K(-\infty)
    \overset{\eqref{eqn-50}}=\frac1{N_p}\sum_{i=1}^{N_p}K_0^i(-\infty).
\end{align}
As $x\rightarrow\infty$, considering the form of $K_0^i(x)$ in \eqref{eqn-59}, we have 
\begin{equation}\label{eqn-43}
    K_0^i(+\infty)=K_0^i(-\infty)+\left(\hat h^{N_p}-C^i\right),
    \end{equation}
since $\textup{erf}(+\infty)=1$. Then,
\begin{align*}
	0\overset{\eqref{eqn-41}}=&p_t^{N_p,\epsilon}K(+\infty)
    \overset{\eqref{eqn-50}}=\frac1{N_p}\sum_{i=1}^{N_p}K_0^i(+\infty)\\
    \overset{\eqref{eqn-43}}=&\frac1{N_p}\sum_{i=1}^{N_p}K_0^i(-\infty)+\frac1{N_p}\sum_{i=1}^{N_p}\left(\hat h^{N_p}-C^i\right)\\
   \overset{\eqref{eqn-73}}=&\frac1{N_p}\sum_{i=1}^{N_p}\left(\hat h^{N_p}-C^i\right),
\end{align*}
which simplifies to
\begin{equation}\label{eqn-75}
	\hat h^{N_p}=\frac1{N_p}\sum_{i=1}^{N_p}C^i.
\end{equation}
Substituting \eqref{eqn-59}-\eqref{eqn-75} back into \eqref{eqn-50}, we arrive at the following result.
\begin{theorem}[Decomposition method for $d=1$]\label{thm-d=1}
	When $d=1$ and $\displaystyle h(x)=\sum_{k=0}^pa_kH_{k}(x)$ is a polynomial of degree $p$, for a given $p\geq1$, where $H_k(x)$ represents the Hermite polynomial of degree $k$. The exact solution to \eqref{eqn-70} is given by
\begin{align}\label{eqn-91}
	K(x)=&\frac1{\sum_{i=1}^{N_p}\mathcal{N}(x;X_t^i,\epsilon)}\\\notag
        &\phantom{a}\cdot\sum_{i=1}^{N_p}\left[\mathcal{N}(x;X_t^i,\epsilon)\sum_{l=0}^{p-1}\hat K^i_lH_l(x)\right.\\\notag
        &\phantom{\sum_{i=1}^{N_p}aa}\left.+\frac{\hat h^{N_p}-C^i}2\textup{erf}\left(\frac{x-X_t^i}{\sqrt{2\epsilon}}\right)\right].
\end{align}
This solution satisfies the boundary condition specified in \eqref{eqn-41}. The coefficients $\{\hat K^i_l\}_{l=0}^{p-1}$ are calculated through backward recursion using \eqref{eqn-92}. The initial values for this recursion are set as $\hat K_{p}^i=\hat K_{p+1}^i=0$. The constant $C^i$ is defined in \eqref{eqn-Ci}.
\end{theorem}

In the one-dimensional case ($d = 1$), the decomposition method under consideration exhibits a computational complexity of $\mathcal{O}(pN_p)$. This indicates that the computational cost increases linearly with both the degree $p$ of the polynomials $h(x)$ involved and the number of particles $N_p$.This linear growth property is a significant advantage when compared to other methods. For instance, both the Galerkin method and the kernel-based method, as reported in \cite{TM:16}, have a computational complexity of $\mathcal{O}(N_p^2)$. A quadratic growth in terms of $N_p$ implies that as the number of particles increases, the computational cost escalates much more rapidly compared to the linear growth of our decomposition method. To further illustrate the superiority of the decomposition method in terms of computational efficiency, numerical comparisons of these methods are presented in Example \ref{ex-4.1}-\ref{ex-4.2}, Section \ref{sec-4.1}.

\subsection{Discussions on the decomposition method in high dimensional scenario}\label{sec-3.c}

The pursuit of exact solutions for \eqref{eqn-3.1} and \eqref{eqn-3.2} in Proposition \ref{prop-1} is an ongoing direction of our research. The Galerkin spectral method, employing multi-variable Hermite polynomials, offers a potential path for solving \eqref{eqn-3.1} to determine the constants $C_j^i$, $j=1,\cdots,d$. But for \eqref{eqn-3.2}, constructing an exact solution that meets the boundary constraint \eqref{eqn-61} is more complex than direct integration in the scalar case, see Section \ref{sec-3.2.2}. This is still among one of the projects we're investigating.

\subsection{Algorithm}\label{AF}

The FPF via decomposition method is detailed in Algorithm \ref{alg-1}.

\begin{algorithm}
\caption{The FPF with decomposition method}
\begin{algorithmic}[1]\label{alg-1}
    \STATE \% Initialization
        \FOR{$i=1$ to $N_p$}
		  \STATE{Sample $X_0^i$ from $p_0^*(x)$}
		\ENDFOR
    \STATE \% The FPF
    \FOR{$n=0$ to $n=T/\Delta t$}
		\STATE Calculate $\hat{h}_t\approx\hat{h}^{N_p}$ in \eqref{eqn-hat hN}
		\FOR{$i=1$ to $N_p$}
            \STATE Generate $N_p$ independent samples $\Delta B^{i}_{t_n}$ from $\N(0,\Delta t)$
            \STATE \% The decomposition method for $K$
            \STATE Calculate $K$ and $u$ for the $i$-th particle $X_{t_n}^{i}$ by \eqref{eqn-91} and \eqref{eqn-u}, respectively
            \STATE Evolve the particles $\{X_{t_n}^{i}\}_{i = 1}^{N_p}$ according to \eqref{eqn-2.20}, i.e.
            \begin{align*}
                X_{t_{n+1}}^{i}=&X_{t_n}^{i}+g(X_{t_n}^i)\Delta t+\sigma(X_{t_n}^i)\Delta B_{t_n}^i\\
                &+K(X_{t_n}^i,t)\Delta Z_{t_k}+u(X_{t_n}^i,t)\Delta t,
            \end{align*} 
        \ENDFOR
    \ENDFOR
\end{algorithmic}
\end{algorithm}

\section{Numeric}\label{sec-4}

\subsection{The comparison of the gain functions}\label{sec-4.1}

In the scalar case, the exact gain function can be obtained by direct integration, if $p_t^*$ is known. That is, 
\begin{equation}\label{eqn-4.1}
 K(x)=-\frac1{p_t^*(x)}\int_{-\infty}^x\left(h(y)-\hat h\right)p_t^*(y)dy.
\end{equation}

We shall first compare the gain functions obtained by various methods with the exact one \eqref{eqn-4.1}.

\begin{example}[Section V, \cite{TMM:17}]\label{ex-1}
Suppose $p_t^*$ is a mixture of two Gaussian distributions, given by $\frac12\N(-\mu,\sigma^2)+\frac12\N(\mu,\sigma^2)$, where $\mu=1$ and $\sigma^2=0.2$. The observation function is defined as $h(x)=x$. Moreover, the number of particles is specified as $N_p=200$.
\end{example}

\begin{figure}[!ht]
    \centering
    \includegraphics[width=0.5\textwidth, trim = 10 200 10 210,clip]{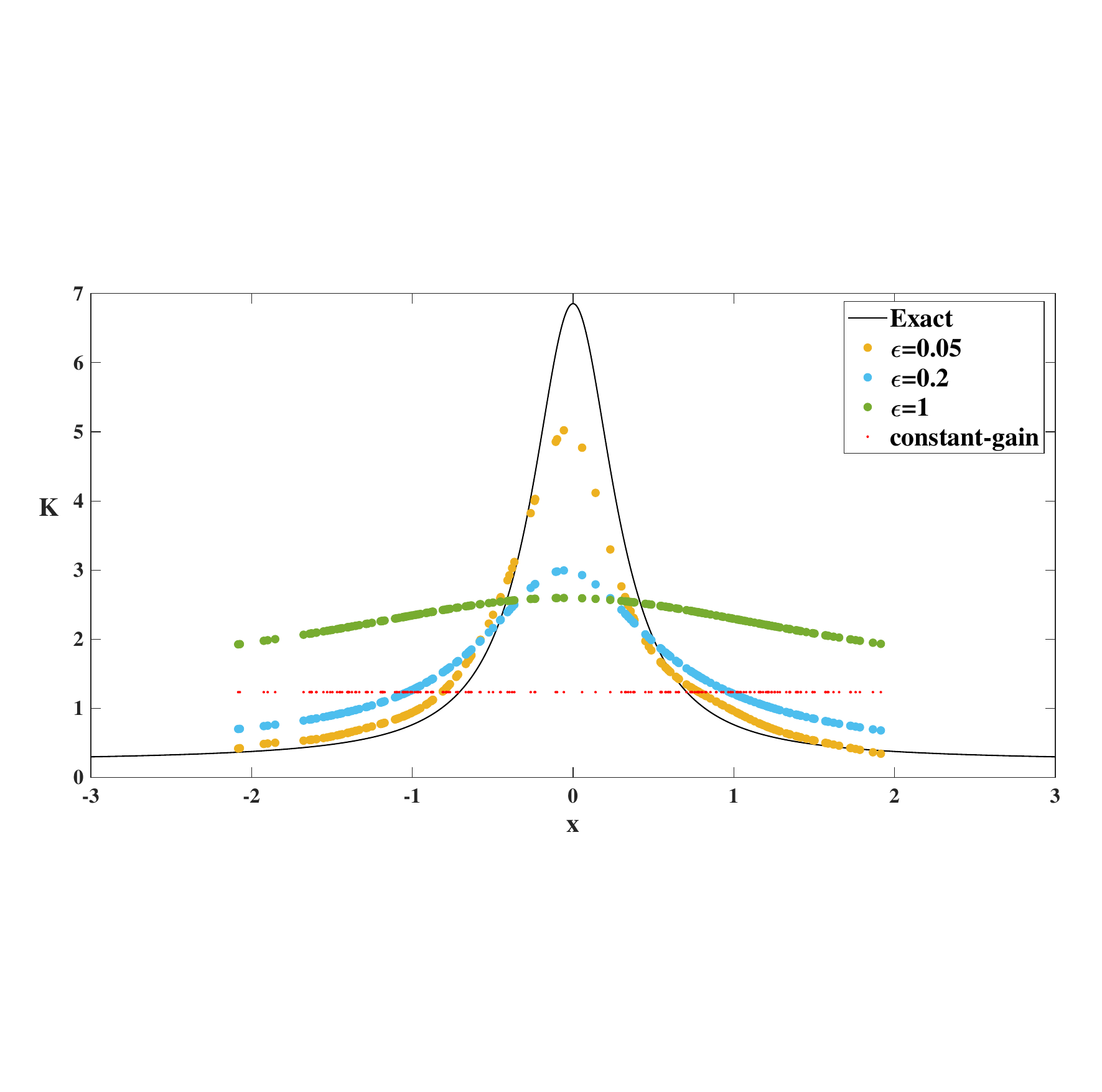}
    \caption{Comparison of the exact gain function and those obtained by the decomposition method with $\epsilon=0.05$, $\epsilon=0.2$, and $\epsilon=1$, respectively.}\label{fig-2}
\end{figure}

In Fig. \ref{fig-2}, we illustrate the gain functions obtained from constant approximation and the decomposition method for $\epsilon=0.05$,  $\epsilon=0.2$ and $\epsilon=1$, respectively. The larger $\epsilon$ is, the flatter the gain function yielded by the decomposition method becomes. In the limit, it converges to a constant value, different from that obtained via the constant-gain approximation. As the value of $\epsilon$ decreases, the approximation derived from the decomposition method becomes more closer to the exact gain function. Nevertheless, when $\epsilon$ is set to be extremely small, the approximation deteriorates significantly because of the presence of variance error.

Next, we shall conduct experiments to verify that the computational time required to obtain the gain function using the decomposition method increases approximately linearly with respect to both the sample size $N_p$ and the degree of the observation function $p$, as we claimed previously in Section \ref{sec-observation}.

\begin{figure}[!ht]
    \centering
    \includegraphics[width=0.5\textwidth, trim = 10 190 10 205,clip]{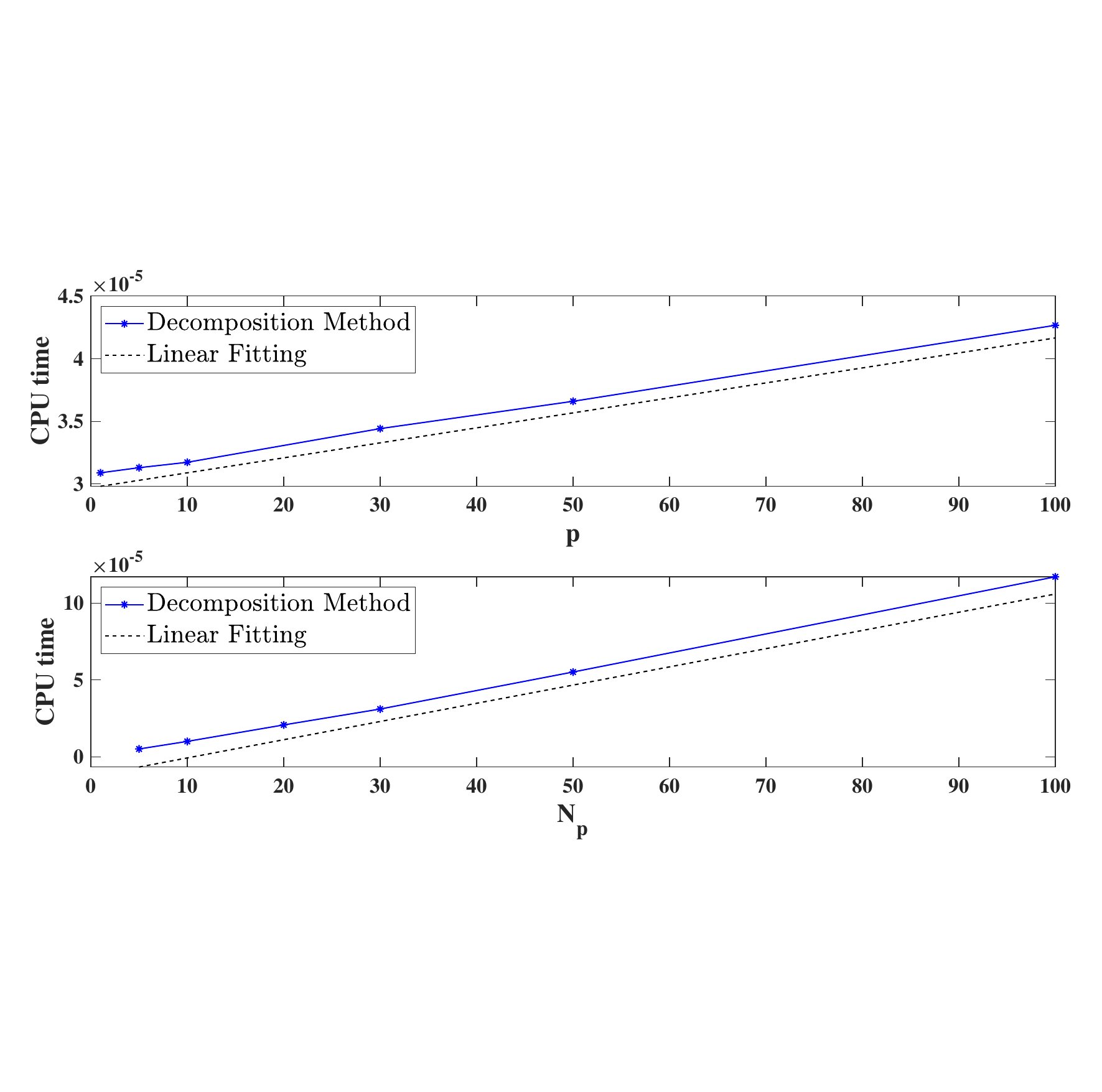}
    \caption{The upper subplot: the relation between CPU times and the degrees of $h(x)$. The lower subplot: the relation between CPU times and the numbers of particles $N_p$.}\label{fig-3&4}
\end{figure}

\begin{example}\label{ex-4.1}
Let the observation functions be $h(x)=H_1(x)$, $H_{5}(x)$, $H_{10}(x)$, $H_{30}(x)$, $H_{50}(x)$, and $H_{100}(x)$. Given that the number of particles is $50$, we shall subsequently utilize the decomposition method to obtain the gain function corresponding to each of these different $h$ functions. Additionally, we shall record the CPU times for computing the gain functions associated with each $h(x)$.
\end{example}

The upper subplot of Fig. \ref{fig-3&4} illustrates the change in the CPU times of the decomposition method as the order of the function $h$ increases. Through the curve fitting, it is readily observable that the CPU times exhibit a linear growth pattern with the increase in the polynomial degree. This observation implies that, when the number of particles is held constant, the computational complexity of the decomposition method is of the order $\mathcal{O}(p)$. Moreover, it is evident that as the polynomial degree rises, the increase in the running time of our algorithm is not substantial. This outcome serves as an evidence for the efficiency of the decomposition method.

\begin{example}\label{ex-4.2}
Let $h(x)=x$ as in Example \ref{ex-1}. The number of particles is $5,10,20,30,50$, and $100$. Then we use the decomposition method to get the gain function for various particle numbers.
\end{example}

In the lower subplot of Fig. \ref{fig-3&4}, as the number of particles increases, the CPU times predominantly exhibit a linear growth tendency. This finding aligns with our prior deductions.

\subsection{A benchmark example}

In this subsection, we numerically solve a NLF problem presented in Section IV.B of \cite{B:18} as a benchmark example. The aim is to demonstrate the tracking capability of our algorithm and to compare its accuracy and efficiency with other methods, such as the FPF with constant-gain approximation and the kernel-based approach.

The benchmark example is modeled as
\begin{equation}\label{NE1}
\left\{\begin{aligned}
    X_{t_{n+1}}=&0.5X_{t_n}+25\frac{X_{t_n}}{1+X_{t_n}^2}+8\cos(1.2n)+B_{t_n}\\
    Z_{t_n}=&0.05X_{t_n}^2+W_{t_n}
\end{aligned}\right.,
\end{equation}
where $B_{t_n}\sim\mathcal{N}(0,10)$, $W_{t_n}\sim\mathcal{N}(0,1)$ and $X_0=0.1$.

The observation function $h(x)$ has coefficients $a_2=\frac1{80}$ and $a_0=\frac1{40}$ as specified in \eqref{eqn-42}. We initiate the backward recursion in \eqref{eqn-92} from $k=3$, with the initial conditions $\hat{K}_{3}^{i}=\hat{K}_{2}^{i}=0$, and then perform the recursion down to $k=0$. This process yields the following results: $\hat{K}_{1}^{i}=\frac{\epsilon}{40}$, $\hat{K}_{0}^{i}=\frac{X_t^i\epsilon}{20}$, and $C^i=\frac{{X_t^i}^2+\epsilon}{20}$. Subsequently, by applying Theorem \ref{thm-d=1}, we are able to obtain the gain function $K(x)$ in an explicit form.

\begin{figure}[!ht]
      \centering
    \includegraphics[width=0.5\textwidth, trim = 10 200 10 220,clip]{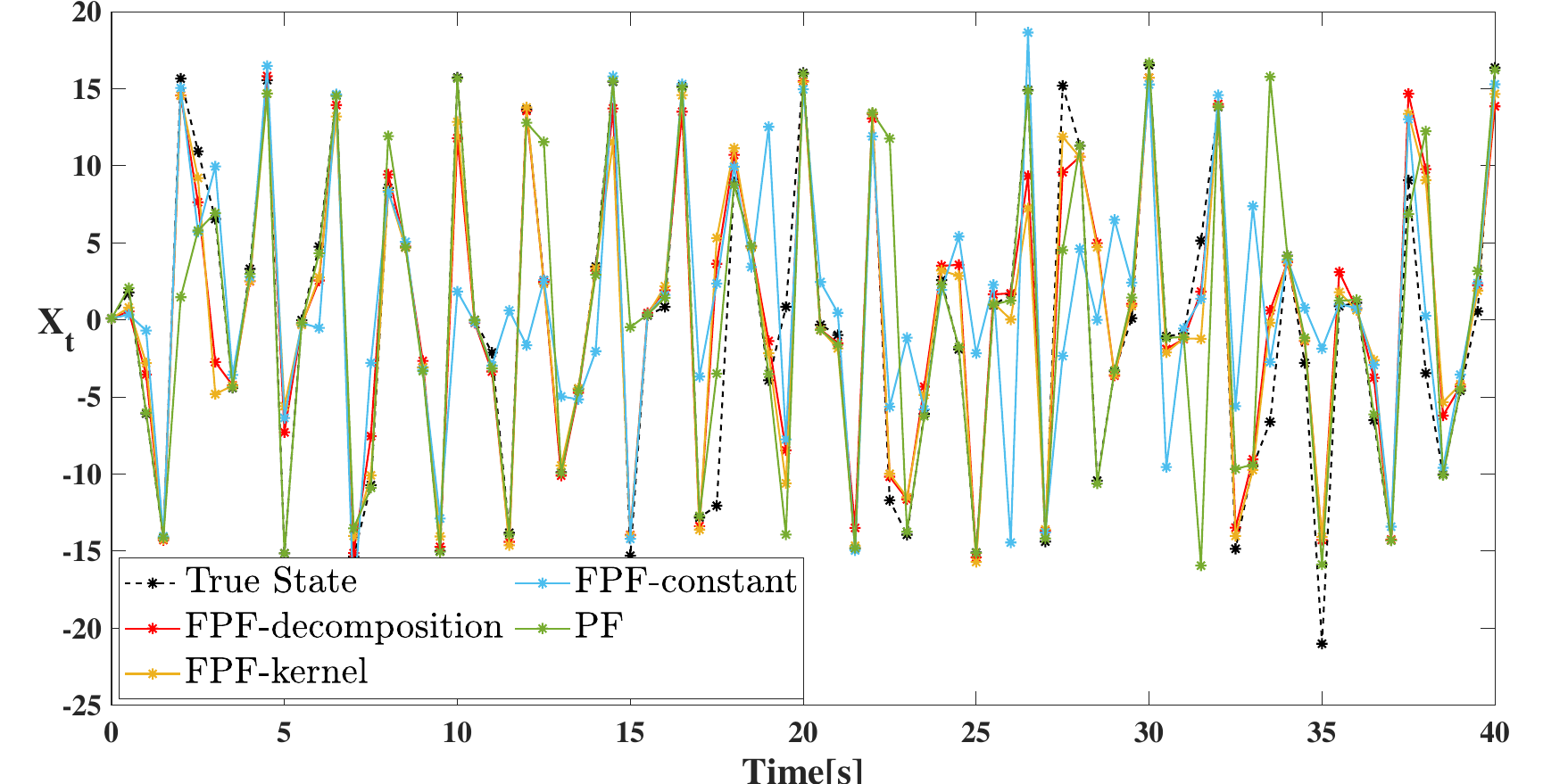}
    \caption{The estimations of the true state (black) obtained by the FPF with the decomposition method (red), the kernel-based approach (yellow), the constant-gain approximation (blue), and the PF (green), respectively.}\label{fig-5}
\end{figure}

\begin{figure}[!ht]
      \centering
      \includegraphics[width=0.5\textwidth, trim = 10 200 10 220,clip]{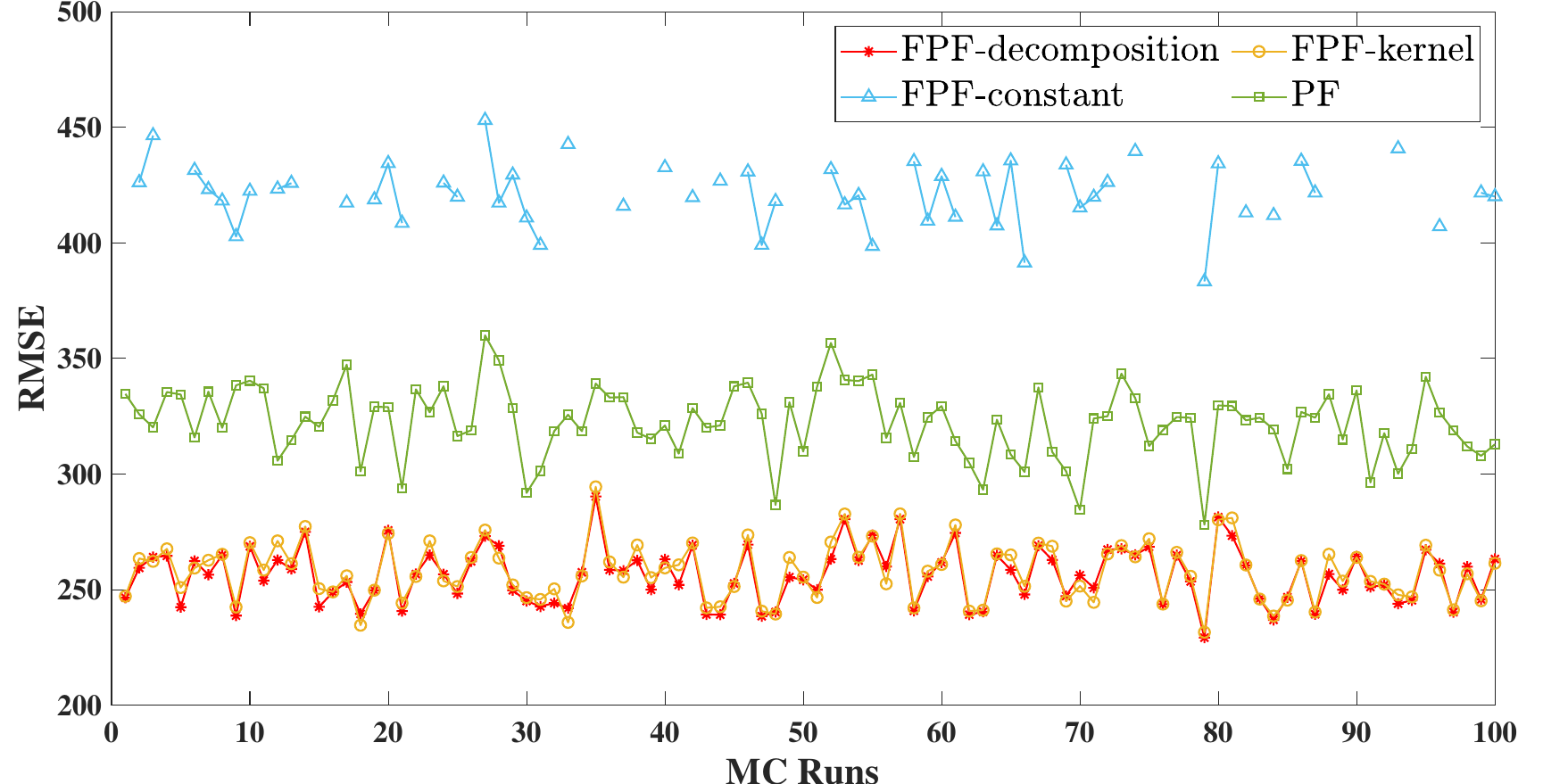}
    \caption{The MSEs averaged over $100$ MC simulations of the FPF with the decomposition method (red), the kernel-based approach (yellow) and the constant-gain approximation (blue), and the PF approximation (green), respectively.}\label{fig-6}
\end{figure}

In this benchmark example, we experiment for the total time $T=40$, and the time step is set to $dt = 0.01$. We generate $50$ particles that follow the standard normal distribution $\N(0,1)$. The true state is generated via the Euler-Maruyama method based on the first SDE in \eqref{NE1}. The parameter $\epsilon$ is set to $0.01$. We first illustrate the performance of the FPF with decomposition method for a single realization. To showcase the effectiveness of this algorithm, we compare the decomposition method with the PF, the FPF with the constant-gain approximation and the kernel-based approach, as depicted in Fig. \ref{fig-5}. To assess the robustness of the algorithm, we perform $100$ Monte-Carlo (MC) simulations, and the results are presented in Fig. \ref{fig-6}. To quantify the error of the algorithm, we utilize the mean square error (MSE). For a single realization, the MSE is defined as 
$\left[\sum\limits_{i=1}^{N_p}(X_{t_n}^i-\hat X_{t_n}^i)^2\right]^\frac{1}{2}$, where $X_{t_n}^i$ and $\hat X_{t_n}^i$ are the true state and the estimation at time $t_n$, respectively. The MSE of decomposition method averaged over $100$ MC simulations is approximately $256.24$, while that of kernel-based approach and the PF with the same number of particles is approximately $257.72$ and $322.36$, respectively. In these $100$ MC simulations, the performance of constant-gain approximation is unstable. Specifically, only $55$ out of the $100$ experiments managed to track the true state. For these $55$ successful experiments, the averaged MSE is roughly $421.49$. The decomposition method, as well as the kernel-based approach, achieves nearly $39\%$ and $20\%$ reductions in MSE compared to the constant-gain approximation and the PF with the same number of particles, respectively. Regarding the computational efficiency, the CPU time averaged over $100$ MC simulations of decomposition method is $0.48$s, that of the kernel-based approach is $0.85$s, of the constant-gain approximation is $0.03$s, and that of the PF is $1.13$s. The decomposition method achieves a $43.53\%$ and a $57.52\%$ reduction in CPU times compared to the kernel-based approach and the PF, respectively.

\section{CONCLUSIONS}

In this paper, we present a novel decomposition method for obtaining the gain function in the FPF. The core concept underlying this method is to decompose the solution of Poisson's equation into two equations that are exactly solvable. Moreover, we conduct a comprehensive comparison of our method with the PF and the FPF with the kernel-based approach and the constant-gain approximation. The theoretical computational complexity of decomposition method is $\mathcal{O}(pN_p)$, which is significantly reduced compared to the kernel-based approach of $\mathcal{O}(N_p^2)$. This is clearly reflected in the computational time in the numerical experiments.

Although the present framework has only been validated in one-dimensional scenarios, it has already yielded favorable experimental results. Its extension to multidimensional systems remains one of our ongoing projects. Future research efforts will be focused on deriving exact solutions for \eqref{eqn-3.1}-\eqref{eqn-3.2} in Proposition \ref{prop-1}. We expect more savings in the CPU times in the high-dimensional NLF problems. This may open up a way to alleviate the curse of dimensionality.

\addtolength{\textheight}{-12cm}   % This command serves to balance the column lengths
                                  % on the last page of the document manually. It shortens
                                  % the textheight of the last page by a suitable amount.
                                  % This command does not take effect until the next page
                                  % so it should come on the page before the last. Make
                                  % sure that you do not shorten the textheight too much.

%%%%%%%%%%%%%%%%%%%%%%%%%%%%%%%%%%%%%%%%%%%%%%%%%%%%%%%%%%%%%%%%%%%%%%%%%%%%%%%%

%%%%%%%%%%%%%%%%%%%%%%%%%%%%%%%%%%%%%%%%%%%%%%%%%%%%%%%%%%%%%%%%%%%%%%%%%%%%%%%%

%%%%%%%%%%%%%%%%%%%%%%%%%%%%%%%%%%%%%%%%%%%%%%%%%%%%%%%%%%%%%%%%%%%%%%%%%%%%%%%%

%\section*{ACKNOWLEDGMENT}

%%%%%%%%%%%%%%%%%%%%%%%%%%%%%%%%%%%%%%%%%%%%%%%%%%%%%%%%%%%%%%%%%%%%%%%%%%%%%%%%

%\bibliographystyle{plain}
%\bibliography{references}

\end{document}